\documentclass[12pt]{article}
\usepackage[cp1251]{inputenc}
\usepackage[english]{babel}
\usepackage[intlimits]{amsmath}
\usepackage{amsfonts}
\usepackage{amssymb}
\usepackage{graphicx}
\usepackage{amsthm}
\usepackage {amsmath}

\begin{document}

\newtheorem{theorem}{Theorem}
\newtheorem{lemma}{Lemma}
\newtheorem{Prop}{Proposition}
\newtheorem{zam}{Remark}
\newtheorem{sled}{Corollary}
\newtheorem{opr}{Definition}
\centerline{\large \bf CONTINUED FRACTIONS}
\centerline{\large \bf WITH ODD PARTIAL QUOTIENTS} 
\vspace*{4mm} \centerline{\Large  Zhabitskaya E.\,N.
\footnote{ Research is supported by RFBR grant No 09-01-00371-a\\
K e y w o r d s: continued fractions, "odd" Euclidean algorithm, distribution function. } }
\vspace*{4mm} ABSTRACT. Consider the representation of a rational number as
a continued fraction, associated with "odd" Euclidean algorithm. In this paper we prove
certain properties for the limit distribution function for sequences of
rationals with bounded sum of partial quotients.

\section{Introduction and main results}
The classical
Euclidean algorithm leads to ordinary continued fraction expansion of a real number
\begin{equation}\label{cf}
x = [ b_0; b_1, b_2, \ldots, b_l, \ldots] =
b_0+\cfrac{1}{b_1+\cfrac{1}{b_{2}
+ \ldots + \cfrac{1}{b_l + \ldots}}},
\end{equation}
where $b_0 \in \mathbb{Z}$, $b_j \in \mathbb{N}$ for $j \geqslant 1.$ 
For rational $x$ this representation is finite.

There are different kinds of Euclidean algorithms (for example,
"by-excess", "centered", "odd" Euclidean algorithms).
Each of them is associated with a kind of continued fraction expansion of a real number
(such fractions can be found in the book \cite{Perron} by O.\,Perron).

Cases of "by-excess" and "centered" Euclidean algorithms were considered in papers 
\cite{zhabitskaya1}, \cite{zhabitskaya2} correspondingly.
In this paper we consider "odd" Euclidean algorithm. This algorithm uses
"odd" division, i.\,e.
\begin{equation*}
a = b q + r, \quad q = 2 \left\lceil \frac{a}{2b} \right\rceil +1, \quad -b < r \leqslant b,
\end{equation*}
and leads to the following representation of a number $x \in \mathbb{Q} \cap [0,1]$:
\begin{equation}\label{frac_odd}
x =\left[1; \frac{\varepsilon_1}{a_1},\ldots,
\frac{\varepsilon_l}{a_l} \right]= 1 +
\cfrac{\varepsilon_1}{a_1 + \cfrac{\varepsilon_2}{a_2+ \ldots +
\cfrac{\varepsilon_l}{a_l}}},
\end{equation}
where all $a_i$ are odd, $\varepsilon_i = \pm 1$ ($\varepsilon_1 = - 1$)
and $a_j + \varepsilon_{j+1} \geqslant 2$  for $j \geqslant
1.$ If the last partial quotient is $a_l=1$,
then $\varepsilon_l = 1$ for uniqueness of the representation.
For irrational $x$ representation (\ref{frac_odd}) is infinite.
We'll call this representation \textit{odd continued fraction}. One can find these fractions
in paper \cite{Rieger}. 
                                                                   
Let us call $F(x)$ by limit distributional function of sequence $\mathcal{M}_n$,
where $\mathcal{M}_n$ is a final subset of segment $[0,1]$, 
if 
$$\lim_{n \to \infty} F_n(x)=F(x),$$
where
$$F_n(x) = \frac{\sharp \left\{ \xi \in \mathcal{M}_n: 
\xi \leqslant x\right\}}{\sharp \mathcal{M}_n}, \quad x \in [0,1].$$ 

We denote the sum of all partial quotients of 
representation (\ref{frac_odd}) of a rational number $x \in [0,1]$ by
$$
S(x) = \sum_{j = 1}^{l} a_j,
$$
and put
$$
\mathcal{M}_n = \left\{x \in \mathbb{Q} \cap [0,1]: S(x) \leqslant n+1\right\}.
$$

In such a way limit distributional function can be defined for any kind of 
continued fraction representation.
For ordinary continued fractions function $F(x)$ coincides with famous
Minkowski's question mark function $?(x)$ (properties of $?(x)$ were
investigated in \cite{Denjoy}, \cite{Salem}).
For regular reduced continued fractions ("by-excess" Euclidean algorithm) and for continued 
fractions with minimal reminders ("centered" Euclidean algorithm) functions $F(x)$ were 
described in papers \cite{zhabitskaya1,zhabitskaya2} correspondingly.
In present paper we consider the function $F(x)$ for odd continued fractions.

The main result is the following theorem.
\begin{theorem}\label{th1}
Suppose that $x \in [0,1]$ is represented in the form (\ref{frac_odd}), 
then
\begin{equation}\label{F}
F(x) = 1-\sum_{i = 1}^{\infty} \frac{E_i}{\lambda^{A_i}},
\end{equation}
where $$E_i = \prod_{ j = 1}^{i}(-\varepsilon_j), \quad
A_i = \sum\limits_{j = 1}^{i}a_j - 1,$$
and $\lambda$ is the unique real root of the equation
\begin{equation*}
\lambda^3 - \lambda^2 - \lambda - 1 = 0.
\end{equation*}
For rational $x$ the sum in formula (\ref{F}) is finite.
\end{theorem}

As a consequence of Theorem \ref{th1} we prove a formula for $F(x)$ in terms of partial quotients of 
ordinary continued fraction.
\begin{sled}\label{sled1}
Suppose that $x \in [0,1]$ is represented in the form (\ref{cf}),
then we have
\begin{equation}\label{F_sled}
F(x) = 1-\sum_{i = 1}^{\infty} (-1)^{i+1}
\frac{c_i}{\lambda^{ S_i(x)-1}},
\end{equation}
where 
\begin{equation*}
c_i = \begin{cases}
 1, \quad b_i \text{~--- odd},\\
 1 + 1/\lambda, \quad b_i \text{~--- even},
\end{cases}
\end{equation*}
$$
S_i(x) = b_1 + \ldots + b_i + \sharp\left\{ j \leqslant i : b_j \text{~--- even} \right\}.
$$
For rational $x$ the sum in formula (\ref{F_sled}) is finite.
\end{sled}

In this paper we also prove the following result.
\begin{Prop}\label{equ_F}
For $x \in [0, 1]$, $n \in \mathbb{N}$ function $F(x)$ satisfies following
functional equations
\begin{gather*}
\frac{1 - F(1 - x)}{\lambda^{2n - 1}} = \frac{1}{\lambda^{2n-2}} - 1 +
F\left(1 - \frac{1}{2n -1 + x}\right),\\
\frac{1 - F(1 - x)}{\lambda^{2n}} = 1 - F\left(1 - \frac{1}{2n + \frac{1}{x}}\right).
\end{gather*}
\end{Prop}

It is more convenient for us to consider representation of $x \in \mathbb{Q} \cap [0,1]$
with the first partial quotient $a_0 = 0$ instead of 1:
\begin{equation}\label{frac_odd_0}
\left[0; \frac{\varepsilon_1}{a_1},\ldots,
\frac{\varepsilon_l}{a_l} \right]:= 
1 - \left[1; \frac{-\varepsilon_1}{a_1},\ldots,
\frac{\varepsilon_l}{a_l} \right]  =
\cfrac{\varepsilon_1}{a_1 + \cfrac{\varepsilon_2}{a_2+ \ldots +
\cfrac{\varepsilon_l}{a_l}}},
\end{equation}
where $\left[1; \frac{-\varepsilon_1}{a_1},\ldots,
\frac{\varepsilon_l}{a_l} \right]$ is representation of number $1 - x$ in the form 
(\ref{frac_odd}).

The limit distributional function corresponding to this representation we denote by $F^{0}(x)$.
For function $F^{0}(x)$ we prove the following results.
\begin{theorem}\label{th_F^0}
Suppose that $x \in [0,1]$ is represented in the form (\ref{frac_odd_0}),
then
\begin{equation}\label{F^0}
F^0(x) = -\sum_{i = 1}^{\infty}\frac{E_i}{\lambda^{A_i}},
\end{equation}
where $$E_i = \prod_{ j = 1}^{i}(-\varepsilon_j), \quad
A_i = \sum\limits_{j = 1}^{i}a_j - 1,$$
and $\lambda$ is the unique real root of the equation
\begin{equation*}
\lambda^3 - \lambda^2 - \lambda - 1 = 0.
\end{equation*}
For rational $x$ the sum in formula (\ref{F^0}) is finite.
\end{theorem}

\begin{Prop}\label{equ_F^0_new}
For $x \in [0, 1]$, $n \in \mathbb{N}$ function $F^0(x)$ satisfies following
functional equations
\begin{gather*}
\frac{F^0(x)}{\lambda^{2n - 1}} = \frac{1}{\lambda^{2n-2}} - F^0\left(\frac{1}{2n -1 + x}\right),\\
\frac{F^0(x)}{\lambda^{2n}} = F^0\left(\frac{1}{2n + \frac{1}{x}}\right).
\end{gather*}
\end{Prop}

\begin{Prop}\label{connect}
For all $x \in [0,1]$ we have 
\begin{equation*}
F(x) = 1 - F^0(1 - x).
\end{equation*}
\end{Prop}

Thus, Theorem \ref{th1} and Proposition \ref{equ_F} follow immediately from 
Theorem \ref{th_F^0}, Proposition \ref{equ_F^0_new} and Proposition \ref{connect}. So
our main aim is to prove results for function $F^0(x)$.

In the end of the paper we prove the following theorem.
\begin{theorem}\label{sing}
Let for $x \in[0,1]$ the derivative $F'(x)$ (finite or infinite)
exists. Then either $F'(x) = 0$ or $F'(x) = \infty$.
\end{theorem}
As function $F(x)$ is monotonic, then by Lebesgue's theorem, the derivative
$F'(x)$ exists and is finite
almost everywhere (in the sense of Lebesgue measure).
That is why $F'(x) = 0$ almost everywhere.
In other words, $F(x)$ is a singular function.

\section{Auxiliary results}
Let us denote by $S^0(x)$ sum of partial quotients of representation (\ref{frac_odd_0})
of a number $x \in \mathbb{Q} \cap [0,1]$.
We define sequences of sets $\mathcal{Y}_n$ and $\mathcal{X}_n$ in the following way:
$$\mathcal{Y}_n :=\left\{x \in \mathbb{Q} \cap [0,1]: S^0(x) \leqslant n+1\right\},$$  
$$\mathcal{X}_k = \left\{x \in \mathbb{Q} \cap [0,1]: S^0(x) =
k+1\right\},$$ 
where $n,k \geqslant 1.$

It is clear that
$$\mathcal{Y}_n = \mathop{\cup}\limits_{1 \leqslant k \leqslant n} \mathcal{X}_k.$$
Suppose that the elements of $\mathcal{Y}_k$ are
arranged in increasing order. The number of elements of
$\mathcal{Y}_n$, $\mathcal{X}_n$ we denote by $Y_n$, $X_n$
correspondingly.

Particularly, $\mathcal{X}_1 =\left\{\frac 1 2\right\}$,
$\mathcal{X}_2 =\left\{\frac 1 3, \frac 2 3\right\}$, $\mathcal{X}_3
=\left\{\frac 1 4, \frac 3 5, \frac 3 4\right\}$,
$\mathcal{X}_4 =\left\{\frac 1 5, \frac 2 7, \frac 2 5, \frac
4 7, \frac 5 8, \frac 4 5 \right\}$. So $X_1 =1$, $X_2 = 2$, $X_3 = 3$, $X_4 =
6$.

\begin{lemma}\label{X_n}
For $n \geqslant 1$ we have
$$X_{n+3} = X_{n+2} + X_{n+1} + X_{n}.$$
\end{lemma}
{\bf Proof.} We construct one-to-one correspondence $\Phi$ between
elements of sets $\mathcal{X}_{n+2} \cup \mathcal{X}_{n+1} \cup
\mathcal{X}_{n}$ and $\mathcal{X}_{n+3}$.

Let $x \in \mathcal{X}_{n+2} \cup \mathcal{X}_{n+1} \cup \mathcal{X}_{n}$,
$x = \left[0; \cfrac{\varepsilon_1}{a_1},\ldots, \cfrac{\varepsilon_l}{a_l} \right]$,
we define $\Phi(x):\mathcal{X}_{n+2} \cup \mathcal{X}_{n+1} \cup \mathcal{X}_{n}
\to \mathcal{X}_{n+3}$ in the following way:
\begin{itemize}
\item In case $x \in \mathcal{X}_{n+2}$
if $a_l = 1$, then
$$\Phi(x)=\left[0; \cfrac{\varepsilon_1}{a_1},
\ldots, \cfrac{\varepsilon_{l-2}}{a_{l-2}}, \cfrac{\varepsilon_{l-1}}{a_{l-1}+2} 
\right] \in \mathcal{X}_{n+3},$$
else
$$\Phi(x)=\left[0; \cfrac{\varepsilon_1}{a_1},
\ldots, \cfrac{\varepsilon_l}{a_l}, \cfrac{1}{1} 
\right] \in \mathcal{X}_{n+3}.$$
\item In case $x \in \mathcal{X}_{n+1}$ if $a_l = 1$, then $\varepsilon_l = 1$ and
$$\Phi(x)=\left[0; \cfrac{\varepsilon_1}{a_1},
\ldots, \cfrac{\varepsilon_{l-1}}{a_{l-1}}, \cfrac{1}{1}, \cfrac{1}{1}, 
\cfrac{1}{1} \right] \in \mathcal{X}_{n+3},$$ else
$$\Phi(x)=\left[0; \cfrac{\varepsilon_1}{a_1},
\ldots, \cfrac{\varepsilon_l}{a_l}, \cfrac{-1}{1}, 
\cfrac{1}{1} \right] \in \mathcal{X}_{n+3}.$$
\item In case $x \in \mathcal{X}_{n}$ if $a_{l} =1$, then
$$\Phi(x)=\left[0; \cfrac{\varepsilon_1}{a_1},
\ldots, \cfrac{\varepsilon_{l-1}}{a_{l-1} + 2}, \cfrac{1}{1}, \cfrac{1}{1} \right] \in \mathcal{X}_{n+3},$$
else
$$\Phi(x)=\left[0; \cfrac{\varepsilon_1}{a_1},
\ldots, \cfrac{\varepsilon_l}{a_l}, \cfrac{-1}{1}, \cfrac{1}{1}, \cfrac{1}{1} \right] \in \mathcal{X}_{n+3}.$$
\end{itemize}

The correspondence $\Phi(x)$ is injective by the construction. Let
us show that it is surjective. For any $y \in \mathcal{X}_{n+3}$,
$y=\left[0; \cfrac{\varepsilon_1}{a_1}, \ldots,
\cfrac{\varepsilon_l}{a_l} \right]$ we find the preimage $x$ of
$y$.

\begin{itemize}
\item If $a_l > 1$ then
$$x = \left[0; \cfrac{\varepsilon_1}{a_1},
\ldots, \cfrac{\varepsilon_{l}}{a_{l}-2},\cfrac{1}{1} \right] \in \mathcal{X}_{n+2}.$$
\item If $a_l = 1$ and $a_{l-1} > 1$,
then $$x = \left[0; \cfrac{\varepsilon_1}{a_1},\ldots,
\cfrac{\varepsilon_{l-1}}{a_{l-1}}\right] \in \mathcal{X}_{n+2}.$$
\item If $a_l = a_{l-1} = 1$, $\varepsilon_{l-1} =-1$, then $a_{l-2} > 1$,
therefore $$x = \left[0; \cfrac{\varepsilon_1}{a_1},\ldots,
\cfrac{\varepsilon_{l-2}}{a_{l-2}}\right] \in \mathcal{X}_{n+1}.$$
\item 
If $a_l = a_{l-1} = 1$, 
$\varepsilon_{l-1} =1$, then 
either $a_{l-2} > 1$ and
$$x = \left[0; \cfrac{\varepsilon_1}{a_1},
\ldots, \cfrac{\varepsilon_{l-2}}{a_{l-2}-2}, \frac{1}{1} \right] \in \mathcal{X}_{n}$$
or
$a_{l-2} = 1$. In this case either $\varepsilon_{l-2} = 1$ and
$$x = \left[0; \cfrac{\varepsilon_1}{a_1},
\ldots, \cfrac{\varepsilon_{l-3}}{a_{l-3}}, \cfrac{1}{1} \right] \in \mathcal{X}_{n+1}$$
or $\varepsilon_{l-2} = -1$, then $a_{l-3}>1$ and
$$x = \left[0; \cfrac{\varepsilon_1}{a_1},
\ldots, \cfrac{\varepsilon_{l-3}}{a_{l-3}} \right] \in \mathcal{X}_{n}$$
\end{itemize}

Lemma is proved.
 \hfill $\blacksquare$

\begin{lemma}
For $n \geqslant 1$ we have
\begin{equation}\label{z_n}
Y_{n+3} = Y_{n+2} + Y_{n+1} + Y_{n} + 2.
\end{equation}
\end{lemma}
{\bf Proof.} By the definition of $\mathcal{Y}_n$ and Lemma \ref{X_n}, we
get
\begin{multline*}
Y_{n+2} + Y_{n+1} + Y_{n} = \\ =
\left(X_{1} + \ldots + X_{n+2}\right)  + \left(X_{1} + \ldots + X_{n+1}\right) +
\left(X_{1} + \ldots + X_{n} \right)= \\ = X_{1} + X_{2} + X_3 + X_{4} + \ldots +
X_{n+3} + \left(X_{1} - X_{3}\right)=  Y_{n+3} - 2.$$
\end{multline*}
\hfill $\blacksquare$

We remind the definition of the Stern-Brocot sequences
$\mathcal{F}_n$, $n = 0,1,2,\ldots$. Consider two-point set
$\displaystyle \mathcal{F}_0 = \left\{ \frac 0 1 , \frac 1 1
\right\}$. Let $n \geqslant 0$ and
$$\mathcal{F}_n =
\left\{0 = x_{0,n} < x_{1,n}< \ldots < x_{N(n),n} = 1 \right\},$$
where $x_{j,n} = p_{j,n}/q_{j,n}$, $(p_{j,n},q_{j,n}) = 1$,
$j = 0,\ldots, N(n)$ and $N(n) = 2^n+1$.
Then
$$\mathcal{F}_{n+1} =\mathcal{F}_{n} \cup Q_{n+1}$$
with
$$Q_{n+1} = \left\{x_{j-1,n}\oplus x_{j,n},\quad j = 1, \ldots, N(n)\right\}.$$
Here
$$\frac{a}{b} \oplus \frac{c}{d} = \frac{a+b}{c+d}$$
is the mediant of fractions $\displaystyle\frac{a}{b}$ and
$\displaystyle\frac{c}{d}$.

It is convenient to represent sequences $\mathcal{F}_n$ by means
of the binary tree $\mathcal{D}^{[0]}$~(Figure~1). This tree is
called Stern-Brocot's tree.
\begin{figure}[h!]
\centerline{\includegraphics{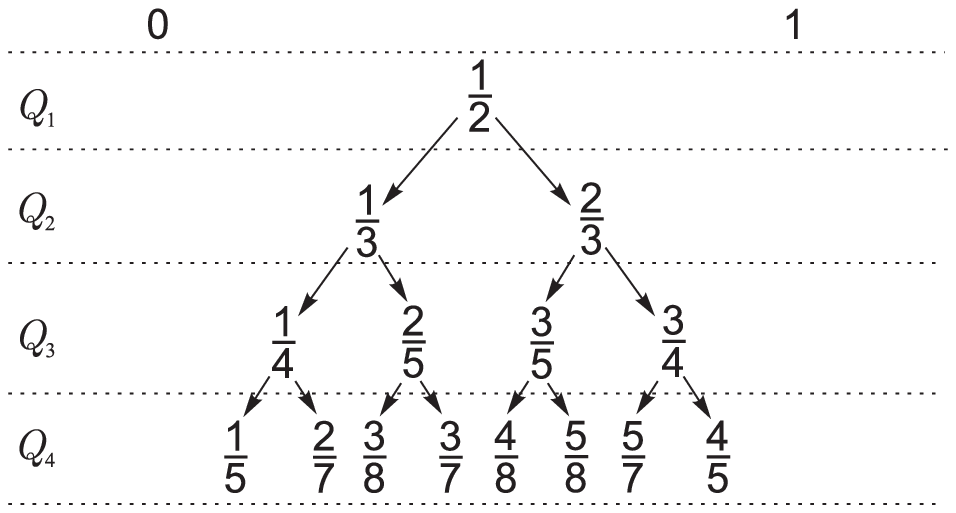}}
\caption{}
\end{figure}
Nodes of the tree are labeled by rationals from interval $(0,1)$ and partitioned into
levels by the following rule:
$n$-th level consists of nodes, labeled by numbers
from ${Q}_n$.

It is possible to distribute nodes of the tree into levels by
another way. For example, we can use such a rule: $n$-th level
consists of nodes labeled by numbers $x$, such that
$S^0(x) = n+1$. 
\begin{figure}[h!]
\centerline{\includegraphics{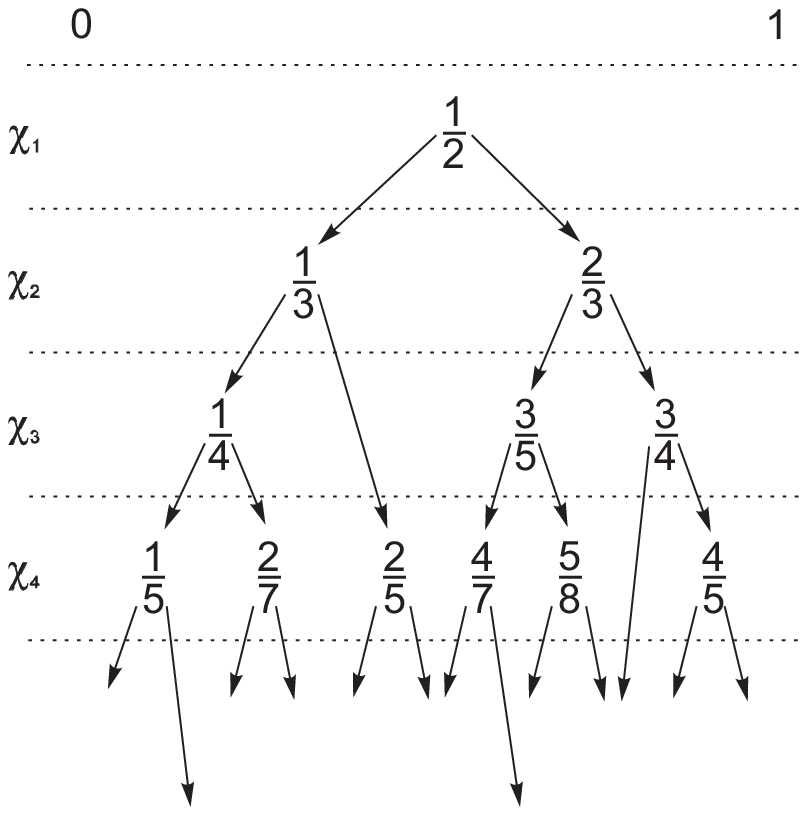}}
\caption{}
\end{figure}
We denote this tree by $\mathcal{D}$~(Figure~2).

{\bf Example.}
\begin{equation*}
\frac{1}{2} = \left[0; \frac{1}{1}, \frac{1}{1} \right],\quad
\frac{1}{3} = \left[0; \frac{1}{3}, \frac{1}{1} \right],\quad
\frac{3}{5} = \left[0; \frac{1}{1}, \frac{1}{1}, \frac{1}{1}, \frac{1}{1} \right],\quad
\frac{2}{5} = \left[0; \frac{1}{3}, \frac{-1}{1}, \frac{1}{1} \right].
\end{equation*}

Any node $\xi$ of the tree $\mathcal{D}$ is a root of a subtree,
which we denote by $\mathcal{D}^{(\xi)}$ ($\mathcal{D} = \mathcal{D}^{(1/2)}$). 
We denote by $D^{(\xi)}_n$ the number
of nodes of $\mathcal{D}^{(\xi)}$ from the level 1 to the level $n$ 
(particularly, $D^{(1/2)}_n = \sharp \mathcal{D}_n$).

Let us consider more detailed structure of the tree $\mathcal{D}$.
From every node $\xi$ of $\mathcal{D}$ we issue two arrows: the
left one and the right one. The left one goes to the node labeled
by ${\xi}^l$ and the right one goes to node labeled by ${\xi}^r$. Note
that if $\xi = x \oplus y$, where $x$, $y$ are consecutive
elements of $\mathcal{F}_n$, then $\xi^l = x \oplus \xi$, $\xi^r =
\xi \oplus y$ (let us call $\xi^l$ and $\xi^r$ successor of $\xi$). 
There are arrows of two kinds: short and long.
Short arrow from $\xi$ to $\eta$, where $\eta \in \{\xi^l, \xi^r\}$, means that 
\begin{equation*}
S^{0}(\eta) -S^{0}(\xi) = 1,
\end{equation*}
and long arrow means that
\begin{equation}\label{S^0 = 2}
S^{0}(\eta) -S^{0}(\xi) = 2.
\end{equation}
Let us call a node with two short arrows by node of first type,
and a node with one short and one long arrow by node of second type.
\begin{Prop}\label{zam}
Suppose that $\xi = \left[0; \cfrac{\varepsilon_1}{a_1},
\ldots, \cfrac{\varepsilon_l}{a_{l}} \right]$.
If $a_l = 1$, then $\xi$ is a node of first type and
if $a_l > 1$, then $\xi$ is a node of second type.
\end{Prop}
{\bf Proof.}
If $a_l = 1$, then
$$\xi^l, \xi^r \in \left\{\left[0; \frac{\varepsilon_1}{a_1},
\ldots, \frac{\varepsilon_{l-1}}{a_{l-1}+2} \right],
\left[0; \frac{\varepsilon_1}{a_1},
\ldots, \frac{\varepsilon_{l}}{a_{l}}, \frac{1}{1} \right]\right\}.$$
If $a_l > 1$, then
$$\xi^l, \xi^r \in \left\{\left[0; \frac{\varepsilon_1}{a_1},
\ldots, \frac{\varepsilon_{l}}{a_{l}}, \frac{1}{1} \right],
\left[0; \frac{\varepsilon_1}{a_1},
\ldots, \frac{\varepsilon_{l}}{a_{l}}, \frac{-1}{1}, \frac{1}{1} \right]\right\}.$$
So the equality (\ref{S^0 = 2}) occurs only in the last case.
\hfill $\blacksquare$

From Proposition \ref{zam} and construction of the tree $\mathcal{D}$
we deduce the following statement.
\begin{lemma}\label{lemma4}
Suppose that $\xi = \left[0; \cfrac{\varepsilon_1}{a_1},
\ldots, \cfrac{\varepsilon_l}{a_{l}} \right]$, then
\begin{equation*}
D^{(\xi)}_n =
\begin{cases}
D^{(1/2)}_{n - S^0(\xi) + 2},  \text{if}\quad a_l = 1\\
D^{(1/3)}_{n - S^0(\xi) + 3},  \text{if}\quad a_l > 1.
\end{cases}
\end{equation*}
\end{lemma}

Note that $D^{(1/2)}_n = Y_n$. For brevity we put $D^{(1/3)}_{n+1} = Z_n$.
From construction of the tree $\mathcal{D}$ it is clear that
\begin{equation}\label{Y_NZ_N}
Y_n = Y_{n-1} + Z_{n-1} +1.
\end{equation}
For $Y_n$ we have recurrence formula (\ref{z_n}). Using equality (\ref{Y_NZ_N}) it is easy to
prove a similar formula for $Z_n$:
\begin{equation*}
Z_{n+3} = Z_{n+2} + Z_{n+1} + Z_{n} + 2.
\end{equation*}
In paper \cite{zhabitskaya2} the author obtained the following result.
\begin{Prop}\label{Prop}
Let $\lambda$ be the unique real root of the equation
\begin{equation*}
\lambda^3 - \lambda^2 - \lambda - 1 = 0,
\end{equation*}
Then
$$\lim_{n \to \infty} \frac{Y_n}{Y_{n+1}} =
\lim_{n \to \infty} \frac{Z_n}{Z_{n+1}} = \frac{1}{\lambda}, \frac{Y_n}{Z_n} = \frac{1}{\lambda - 1}.$$
\end{Prop}

\section{Proof of results for function $F^0(x)$}
In order to prove Theorem 2 we need the following lemma.
\begin{lemma}\label{lemma3}
Suppose that $x = \left[0; \cfrac{\varepsilon_1}{a_1},
\ldots, \cfrac{\varepsilon_l}{a_{l}}, \ldots \right]$, 
then
\begin{equation*}
F^0_n(x) = -\left(\sum_{i = 1}^{\infty} E_i D^{\left[0; 
\frac{\varepsilon_1}{a_1},
\ldots, \frac{\varepsilon_i}{a_{i}}, \frac{1}{1} \right]}_n\right) \bigl/ D^{1/2}_n,
\end{equation*}
where $$E_i = \prod_{ j = 1}^{i}(-\varepsilon_j).$$
\end{lemma}
{\bf Proof.}
We will prove the lemma by induction on the length $l$ of odd continued fraction 
representation of
$x = \left[0; \frac{\varepsilon_1}{a_1},
\ldots, \frac{\varepsilon_l}{a_{l}} \right]$. 

For $l = 1$ we have
$$F^0_n([0; 1/a_1]) = \frac{\sharp \left\{ \xi \in Y_n: 
\xi \leqslant 1/a_1\right\}}{ Y_n} = 
\frac{D^{0/1 \oplus 1/a_1}_n}{D^{1/2}_n} = 
\frac{D^{\left[0; \frac{1}{a_1}, \frac{1}{1}\right]}_n}{D^{1/2}_n}.$$

Now suppose that $l = i+1$ and the lemma is true for $l \leqslant i$. 

Put $a = \left[0; \frac{\varepsilon_1}{a_1},
\ldots, \frac{\varepsilon_{i}}{a_{i}}, \frac{\varepsilon_{i+1}}{a_{i+1}} 
\right]$,
$b = \left[0; \frac{\varepsilon_1}{a_1},
\ldots, \frac{\varepsilon_{i}}{a_{i}} \right]$,
$c = \left[0; \frac{\varepsilon_1}{a_1},
\ldots, \frac{\varepsilon_{i-1}}{a_{i-1}} \right]$
(if the last partial quotient of $a$, $b$ or $c$ is $-1/1$, we replace it by $1/1$
and decrease previous partial quotient by 2).
Suppose that $c < b$ (in case $c < b$ the proof is analogously).
By assumption of induction we have
$$F^0_n(b) - F^0_n(c) = \frac{D^{\left[0; \frac{\varepsilon_1}{a_1},
\ldots, \frac{\varepsilon_{i}}{a_{i}}, \frac{1}{1} \right]}_n}{D^{1/2}_n}.$$

Taking into account the fact that $a$ and $b$ are consecutive elements
of an element of Stern-Brocot sequence (as consecutive convergents of number $x$) we have 
$$
a \oplus b = \left[0; \frac{\varepsilon_1}{a_1},
\ldots, \frac{\varepsilon_{i}}{a_{i}}, 
\frac{\varepsilon_{i+1}}{a_{i+1}}, \frac{1}{1} \right].
$$
By definition of $F_n$ 
\begin{multline*}
F^0_n\left(\max{(a,b)}\right) - F^0_n\left(\min{(a,b)}\right)= 
\frac{\sharp \left\{ \xi \in Y_n: 
\min{(a,b)} < \xi \leqslant \max{(a,b)} \right\}}{ Y_n} =\\ = 
\frac{D^{a \oplus b}_n}{D^{1/2}_n} = 
\frac{D^{\left[0; \frac{\varepsilon_1}{a_1},
\ldots, \frac{\varepsilon_{i}}{a_{i}}, 
\frac{\varepsilon_{i+1}}{a_{i+1}}, \frac{1}{1} \right]}_n}{D^{1/2}_n}.
\end{multline*}  

For $\varepsilon_{i+1} = 1$ we have $a \in (c, b)$, so
$$\sharp\{a <\xi \leqslant b\} = 
\sharp\{c <\xi \leqslant b\}-
\sharp\{a <\xi \leqslant b\}$$
and
$$
F^0_n(a) - F^0_n(c) = \frac{D^{\left[0; \frac{\varepsilon_1}{a_1},
\ldots, \frac{\varepsilon_{i}}{a_{i}}, 
  \frac{1}{1} \right]}_n - D^{\left[0; \frac{\varepsilon_1}{a_1},
\ldots, \frac{\varepsilon_{i}}{a_{i}}, 
 1/a_{i+1},  \frac{1}{1} \right]}_n}{D^{1/2}_n}.
$$
For $\varepsilon_{i+1} = -1$ we have $a > b $, so
$$\sharp\{b <\xi \leqslant a\} = 
\sharp\{c <\xi \leqslant b\}+\sharp\{b <\xi \leqslant a\}$$
and
$$
F^0_n(a) - F^0_n(c)= \frac{D^{\left[0; \frac{\varepsilon_1}{a_1},
\ldots, \frac{\varepsilon_{i}}{a_{i}}, 
 1/1 \right]}_n + D^{\left[0; \frac{\varepsilon_1}{a_1},
\ldots, \frac{\varepsilon_{i}}{a_{i}}, 
 -1/a_{i+1}, \frac{1}{1} \right]}_n}{D^{1/2}_n}.
$$
 \hfill $\blacksquare$

{\bf Proof of the Theorem 2.}
By Lemma \ref{lemma4} and Proposition \ref{Prop} we have
$$
\frac{D^{\left[0; {\varepsilon_1}/{a_1},
\ldots, {\varepsilon_i}/{a_{i}}, 1/1 \right]}_n}{D^{1/2}_n} = 
\frac{Y_{n - (\sum_{j = 1}^{i} a_j + 1)+2}}{Y_n} = 
\frac{1}{\lambda^{\sum_{j = 1}^{i} a_j - 1}}.
$$
So formula (\ref{F^0}) follows immediately from Lemma \ref{lemma3}.
 \hfill $\blacksquare$

%
{\bf Proof of Proposition 2.}
Suppose that $x = \left[0; \frac{\varepsilon_1}{a_1}, \ldots, \frac{\varepsilon_n}{a_n}, \ldots \right]$
is representation of $x$ in the form (\ref{frac_odd_0}).
Then we have
$$\frac{1}{2n-1 + x} = 
\left[0; \frac{1}{2n -1}, \frac{\varepsilon_1}{a_1}, \ldots, \frac{\varepsilon_n}{a_n}, \ldots \right]$$
and
$$\frac{1}{2n + \frac{1}{x}} = 
\left[0; \frac{\varepsilon_1}{2n + a_1}, \frac{\varepsilon_2}{ a_2}, \ldots, 
\frac{\varepsilon_n}{a_n}, \ldots \right].$$
Now it is only left to apply Theorem \ref{th1} to these numbers.

{\bf Proof of Proposition 3.}
Suppose that $x \in \mathbb{Q}$ and
$$
x = \left[1;\frac{\varepsilon_1}{a_1}, \ldots, \frac{\varepsilon_n}{a_n}\right]
$$
is representation of $x$ in the form (\ref{frac_odd}).

Representation (\ref{frac_odd_0}) is connected with (\ref{frac_odd}) in the following way
\begin{equation*}
1 - \left[1; \frac{\varepsilon_1}{a_1},\ldots,
\frac{\varepsilon_l}{a_l} \right] =\left[0; \frac{-\varepsilon_1}{a_1},\ldots,
\frac{\varepsilon_l}{a_l} \right].
\end{equation*}
So we have
$$
1 - x = \left[0;\frac{-\varepsilon_1}{a_1}, \ldots, \frac{\varepsilon_n}{a_n}\right].
$$
Consequently 
$$S(x) = 1 + a_1 + \ldots + a_n = 1 + S^0(1 - x).$$
That is why 
we have
$$
F_n(x) = \frac{\sharp \left\{ \xi \in \mathcal{M}_{n}: \xi \leqslant x 
\right\}}{\sharp \left\{ \xi \in \mathcal{M}_n \right\}}=
\frac{\sharp \left\{ 1 - \xi \in \mathcal{Y}_{n-1}: 1 - \xi \geqslant 1 - x 
\right\}}{ Y_{n-1} } = 1 - F^0_n(1 - x).
$$
For irrational $x\in [0,1]$ we should take into account continuity of considered functions.
\hfill $\blacksquare$

\section{Proof of Corollary 1}
Function $F(x)$ can be also expressed in terms of partial quotients of ordinary continued fraction.

Suppose that $x \in \mathbb{Q} \cap [0,1]$ is represented in the form of ordinary 
continued fraction:
$$
x = [0; b_1, \ldots, b_n].
$$
We describe the algorithm for converting this fraction into a fraction of the form (\ref{frac_odd}).
For the fist $i$ such that $b_i$ is even we use one of the following identities:
\begin{gather}  \label{i1}
b_i + \cfrac{1}{b_{i+1} + \alpha} = 
b_i + 1 - \cfrac{1}{1 + \cfrac{1}{b_{i+1} - 1 + \alpha}}, \\
\label{i2}
b_i + \cfrac{1}{1 + \cfrac{1}{b_{i+2} + \alpha}} = 
b_i + 1 - \cfrac{1}{b_{i+2} + 1 + \alpha},
\end{gather} 
where $\alpha$ is the "tail" of the fraction.

In case $b_{i+1} > 1$ we use identity (\ref{i1}) 
whereas in case $b_{i+1} = 1$ we use identity (\ref{i2}). 
Then we apply the same procedure to the obtained fraction.

To prove Corollary 1 we should apply Theorem \ref{th1} to the result of the procedure described 
above. Taking into account the fact that in case $b_i$ is even, if $b_{i+1} > 1$ we have
\begin{multline*}
\frac{1}{\lambda^{b_1 + \ldots + b_{i-1} + (b_i +1) - 1}} +
\frac{1}{\lambda^{b_1 + \ldots + b_{i-1} + (b_i +1) + 1 - 1}} 
-\frac{1}{\lambda^{b_1 + \ldots + b_{i-1} + (b_i +1) + 1 + (b_{i+1} - 1) - 1}} 
=\\=
\frac{1}{\lambda^{b_1 + \ldots + b_{i-1} + b_i + 1 - 1}} \left(1 + \frac{1}{\lambda}\right) - 
\frac{1}{\lambda^{b_1 + \ldots + b_{i-1} + b_i + b_{i + 1} + 1 - 1}}
\end{multline*}
and if $b_{i+1} = 1$ then
\begin{multline*}
\frac{1}{\lambda^{b_1 + \ldots + b_{i-1} + (b_i +1) - 1}} 
+\frac{1}{\lambda^{b_1 + \ldots + b_{i-1} + (b_i +1) + (b_{i+2} +1)- 1}}  
=\\=
\frac{1}{\lambda^{b_1 + \ldots + b_{i-1} + (b_i + 1) - 1}} \left(1 + \frac{1}{\lambda}\right) - 
\frac{1}{\lambda^{b_1 + \ldots + b_{i-1} + (b_{i} + 1 )+ 1 - 1}} +\\
+\frac{1}{\lambda^{b_1 + \ldots + b_{i-1} + (b_i +1) + 1 + b_{i+2} - 1}}.  
\end{multline*}
So we get formula (\ref{F_sled}).
\hfill $\blacksquare$

\section{Singularity of the function $F(x)$}
In this section we prove Theorem \ref{sing}.
At first we consider the  case $x \in \mathbb{Q}$.
\begin{lemma}
For rational $x \in [0,1]$ we have $F'(x) = 0$.
\end{lemma}
{\bf Proof.}
As $x \in \mathbb{Q}$, so there exists such $n$, that
$a/b \in \mathcal{X}_n$.
By $p/q$, $p'/q'$ we denote the left and the right neighbouring to $a/b$ elements in $\mathcal{Y}_n$
correspondingly.
Sequences of mediants
$$
\xi_k = \frac{k a + p}{k b + q}, \quad
{\xi'}_k = \frac{k a + p'}{k b + q'}.
$$
converge to $a/b$ from the left and from the right correspondingly as $k\to \infty$.

For consecutive elements $x$, $y$ of $\mathcal{Z}_n$ the ration
\begin{equation*}
\left(F(x \oplus y) - F(x)\right):\left(F(y) - F(x \oplus y)\right) 
\end{equation*} 
can take values $\left\{\frac{\lambda}{\lambda+1}, \frac{1}{\lambda+1},
\frac{\lambda-1}{\lambda}, \frac{1}{\lambda}\right\}$,
so
\begin{multline*}
0 \leqslant F(x) - F(\xi_k) \leqslant \left(\max \left\{\frac{\lambda}{\lambda+1}, \frac{1}{\lambda+1},
\frac{\lambda-1}{\lambda}, \frac{1}{\lambda}\right\}\right)^k \left(F(x) - F(p/q)\right) = \\=
\left(\frac{\lambda}{\lambda+1}\right)^k \left(F(x) - F(p/q)\right).
\end{multline*}
Analogously
\begin{equation*}
0 \leqslant F({\xi'}_k) - F(x) \leqslant 
\left(\frac{\lambda}{\lambda+1}\right)^k \left(F(p'/q') - F(x)\right).
\end{equation*}
Taking into account the fact that
$$
{\xi}'_k - \xi_k  = \frac{2k + 1}{(k b + q')(k b + q)} \geqslant \frac{2k +1}{k^2(b + p')(b + p)}
$$
we have
\begin{multline*}
0 \leqslant \lim_{k \to infty}\frac{F({\xi}'_k) - F(\xi_k)}{{\xi}'_k - \xi_k} \leqslant
\lim_{k \to infty} \frac{\left(\frac{\lambda}{\lambda+1}\right)^k \left(F(p'/q') - F(p/q)\right)}
{\frac{2k +1}{k^2(b + p')(b + p)}} = 0
\end{multline*}
\hfill $\blacksquare$

Now we should prove the Theorem for irrational $x \in [0,1]$.

Given $n$ 
we can find two consecutive elements $p_n/q_n<p'_n/q'_n$ from the set $\mathcal{F}_n$
such that $p_n/q_n<x<p'_n/q'_n$. In such a way we obtain an
infinite sequence of pairs of elements $\{p_n/q_n, p'_n/q'_n\}$,
converging to $x$ from the left and from the right correspondingly.
Note that $p_n/q_n$, $p'_n/q'_n$  are always among the  intermediate
and convergent fractions to $x$ in the sense of ordinary continued fraction.
\begin{lemma}\label{cases} There are two possibilities
\begin{enumerate}
\item
There are infinitely many $i$ such that for three consecutive pairs 
$\left\{\frac{p_i}{q_i}, \frac{{p_i}'}{{q_i}'}\right\}$,
$\left\{\frac{p_{i+1}}{q_{i+1}}, \frac{{p_{i+1}}'}{{q_{i+1}}'}\right\}$,
$\left\{\frac{p_{i+2}}{q_{i+2}}, \frac{{p_{i+2}}'}{{q_{i+2}}'}\right\}$
the following equalities hold:
\begin{gather}
\frac{F\left(\frac{p_{i+1}'}{q_{i+1}'}\right)-F\left(\frac{p_{i+1}}{q_{i+1}}\right)}  
{F\left(\frac{p_{i}'}{q_{i}'}\right)-F\left(\frac{p_{i}}{q_{i}}\right)} \in 
\left\{\frac{\lambda - 1}{\lambda}, \frac{1}{\lambda}\right\} \label{1}\\ 
\frac{F\left(\frac{p_{i+2}'}{q_{i+2}'}\right)-F\left(\frac{p_{i+2}}{q_{i+2}}\right)}  
{F\left(\frac{p_{i+1}'}{q_{i+1}'}\right)-F\left(\frac{p_{i+1}}{q_{i}+1}\right)} \in 
\left\{\frac{\lambda - 1}{\lambda}, \frac{1}{\lambda}\right\} \label{2} 
\end{gather}
\item
There are infinitely many $i$ such that for three consecutive pairs 
$\left\{\frac{p_i}{q_i}, \frac{{p_i}'}{{q_i}'}\right\}$,
$\left\{\frac{p_{i+1}}{q_{i+1}}, \frac{{p_{i+1}}'}{{q_{i+1}}'}\right\}$,
$\left\{\frac{p_{i+2}}{q_{i+2}}, \frac{{p_{i+2}}'}{{q_{i+2}}'}\right\}$
the following equalities hold:
\begin{gather}
\frac{F\left(\frac{p_{i+1}'}{q_{i+1}'}\right)-F\left(\frac{p_{i+1}}{q_{i+1}}\right)}  
{F\left(\frac{p_{i}'}{q_{i}'}\right)-F\left(\frac{p_{i}}{q_{i}}\right)} \in 
\left\{\frac{\lambda}{\lambda +1}, \frac{1}{\lambda+1}\right\} \label{3}\\ 
\frac{F\left(\frac{p_{i+2}'}{q_{i+2}'}\right)-F\left(\frac{p_{i+2}}{q_{i+2}}\right)}  
{F\left(\frac{p_{i+1}'}{q_{i+1}'}\right)-F\left(\frac{p_{i+1}}{q_{i}+1}\right)} \in 
\left\{\frac{\lambda - 1}{\lambda}, \frac{1}{\lambda}\right\} \label{4} 
\end{gather}
\end{enumerate}
\end{lemma}
{\bf Proof of the Lemma.}
Let $x$, $y$ be consecutive elements of $\mathcal{F}_n$ for some $n$.
Then by Proposition \label{Prop}
\begin{multline*}
\left(F(x \oplus y) - F(x)\right):\left(F(y) - F(x \oplus y)\right) =\\=
\begin{cases} 
\lambda -1 \quad \text{or} \quad \frac{1}{\lambda-1}, \,\text{if} \quad
x \oplus y\text{~--- vertex of the first type,}\\
\lambda \quad \text{or} \quad \frac{1}{\lambda}, \,\text{if} \quad
x \oplus y\text{~--- vertex of the second type.}
\end{cases}
\end{multline*}

Suppose that fist case of the Lemma is not hold true. Then there are infinitely many $i$ such that
(\ref{3}) holds. But if $\frac{p_i}{q_i} \oplus \frac{p'_i}{q'_i}$ is a vertex of second type, 
then both $\left(\frac{p_i}{q_i} \oplus \frac{p'_i}{q'_i}\right)^l$ and
$\left(\frac{p_i}{q_i} \oplus \frac{p'_i}{q'_i}\right)^r$ are vertexes of first type.
That is why (\ref{4}) holds.

{\bf Proof of the Theorem.}
Suppose that $F' (x) = a$, where $a$ is finite and $a\not = 0$.
By definition of derivative we have
\begin{equation*}
\lim_{i \to \infty} \frac{F\left(\frac{p_{i}'}{q_{i}'}\right)-F\left(\frac{p_{i}}{q_{i}}\right)}
{\frac{p_{i}'}{q_{i}'}-\frac{p_{i}}{q_{i}}}  = a \not = 0,
\end{equation*}
consequently
\begin{equation}
\frac{F\left(\frac{p_{i+1}'}{q_{i+1}'}\right)-F\left(\frac{p_{i+1}}{q_{i+1}}\right)}
{F\left(\frac{p'_{i}}{q'_{i}}\right)-F\left(\frac{p_{i}}{q_{i}}\right)} \sim 
\frac{\frac{p_{i+1}'}{q_{i+1}'}-\frac{p_{i+1}}{q_{i+1}}}{\frac{p_{i}'}{q_{i}'}-\frac{p_{i}}{q_{i}}}. 
\end{equation}
And since
$\left\{\frac{p_{i+1}}{q_{i+1}}, \frac{p'_{i+1}}{q'_{i+1}}\right\}$ is either
$\left\{\frac{p_{i}}{q_{i}}, \frac{p_i +p'_{i}}{q_i + q'_{i}}\right\}$ or
$\left\{\frac{p_i +p'_{i}}{q_i + q'_{i}}, \frac{p'_{i}}{q'_{i}}\right\}$ then
\begin{equation*}
\frac{\frac{p_{i+1}'}{q_{i+1}'}-\frac{p_{i+1}}{q_{i+1}}}
{\frac{p_{i}'}{q_{i}'}-\frac{p_{i}}{q_{i}}}= 
\frac{\frac{1}{q_{i+1}'q_{i+1}}}
{\frac{1}{q_{i}'q_{i}}} \in \left\{\frac{{q_i}'}{q_i + q'_i}, \frac{{q_i}}{q_i + q'_i}\right\}.
\end{equation*}
As
$$
\frac{{q_{i_k}}'}{q_{i_k} + q'_{i_k}} + \frac{{q_{i_k}}}{q_{i_k} + q'_{i_k}} = 1
$$
so if the first case of lemma \ref{cases} holds for a sequence $\{i_k\}_{k = 1}^{\infty}$, then either
\begin{equation}\label{e1}
\lim_{{i_k} \to \infty}\frac{{q_{i_k}}'}{q_{i_k} + q'_{i_k}} = \frac{\lambda - 1}{\lambda}\quad \text{and} \quad
\lim_{{i_k} \to \infty}\frac{{q_{i_k}}}{q_{i_k} + q'_{i_k}} = \frac{1}{\lambda}
\end{equation}
or
\begin{equation}\label{e2}
\lim_{{i_k} \to \infty}\frac{{q_{i_k}}'}{q_{i_k} + q'_{i_k}} = \frac{1}{\lambda}\quad \text{and} \quad
\lim_{{i_k} \to \infty}\frac{{q_{i_k}}}{q_{i_k} + q'_{i_k}} = \frac{\lambda - 1}{\lambda}.
\end{equation}
If the second case holds, then either
\begin{equation}\label{e3}
\lim_{{i_k} \to \infty}\frac{{q_{i_k}}'}{q_{i_k} + q'_{i_k}} = \frac{\lambda}{\lambda + 1}\quad \text{and} \quad
\lim_{{i_k} \to \infty}\frac{{q_{i_k}}}{q_{i_k} + q'_{i_k}} = \frac{1}{\lambda + 1}
\end{equation}
or
\begin{equation}\label{e4}
\lim_{{i_k} \to \infty}\frac{{q_{i_k}}'}{q_{i_k} + q'_{i_k}} = \frac{1}{\lambda +1}\quad \text{and} \quad
\lim_{{i_k} \to \infty}\frac{{q_i}}{q_{i_k} + q'_{i_k}} = \frac{\lambda}{\lambda +1}.
\end{equation}
Analogously the pair 
$\left\{\frac{p_{i+2}}{q_{i+2}}, \frac{p'_{i+2}}{q'_{i+2}}\right\}$ can take one of values
$\left\{\frac{p_{i}}{q_{i}}, \frac{2p_{i} + p'_{i}}{2q_{i} + q'_{i}}\right\}$,
$\left\{\frac{2p_{i} + p'_{i}}{2q_{i} + q'_{i}}, \frac{p_{i} + p'_{i}}{q_{i}+ q'_{i}}\right\}$,
$\left\{\frac{p_{i} + p'_{i}}{q_{i}+ q'_{i}}, \frac{p_{i} + 2p'_{i}}{q_{i} + 2q'_{i}}\right\}$
or
$\left\{\frac{p_{i} + 2p'_{i}}{q_{i} + 2q'_{i}}, \frac{p'_{i}}{q'_{i}}\right\}$. So
\begin{equation*}
\frac{\frac{p_{i+2}'}{q_{i+2}'}-\frac{p_{i+2}}{q_{i+2}}}
{\frac{p_{i+1}'}{q_{i+1}'}-\frac{p_{i+1}}{q_{i+1}}}=
\frac{\frac{1}{q_{i+2}'q_{i+2}}}
{\frac{1}{q_{i+1}'q_{i+1}}} =
\begin{cases} 
\displaystyle \frac{{q_i}}{2q_i + q'_i}\quad \text{or}\quad \frac{{q_i + q'_i}}{2q_i + q'_i},
\quad \text{if} \quad
x < \frac{p_{i} + p'_{i}}{q_{i} + q'_{i}}\\ 
\displaystyle \frac{{q_i + q'_i}}{q_i + 2q'_i}\quad \text{or}\quad \frac{{q'_i}}{q_i + 2q'_i},
 \quad\text{if}  \quad
x > \frac{p_{i} + p'_{i}}{q_{i} + q'_{i}}. 
\end{cases}
\end{equation*}
So in the first case of lemma \ref{cases} for $x < \frac{p_{i} + p'_{i}}{q_{i} + q'_{i}}$ we have either
\begin{equation}\label{e5}
\lim_{{i_k} \to \infty}\frac{{q_{i_k}}}{2q_{i_k} + q'_{i_k}} = \frac{\lambda - 1}{\lambda}\quad \text{and}\quad 
\lim_{{i_k} \to \infty}\frac{{q_{i_k} + q'_{i_k}}}{2q_{i_k} + q'_{i_k}} = \frac{1}{\lambda}
\end{equation}
or
\begin{equation}\label{e6}
\lim_{{i_k} \to \infty}\frac{{q_{i_k}}}{2q_{i_k} + q'_{i_k}} = \frac{1}{\lambda}\quad \text{and}\quad 
\lim_{{i_k} \to \infty}\frac{{q_{i_k} + q'_{i_k}}}{2q_{i_k} + q'_{i_k}} = \frac{\lambda - 1}{\lambda},
\end{equation}
and for $x > \frac{p_{i} + p'_{i}}{q_{i} + q'_{i}}$ we have either
\begin{equation}\label{e7}
\lim_{{i_k} \to \infty}\frac{{q_{i_k} + q'_{i_k}}}{q_{i_k} + 2q'_{i_k}} = \frac{\lambda - 1}{\lambda}\quad \text{and}\quad 
\lim_{{i_k} \to \infty}\frac{{q'_{i_k}}}{q_{i_k} + 2q'_{i_k}}= \frac{1}{\lambda}
\end{equation}
or
\begin{equation}\label{e8}
\lim_{{i_k} \to \infty}\frac{{q_{i_k} + q'_{i_k}}}{q_{i_k} + 2q'_{i_k}} = \frac{1}{\lambda}\quad \text{and}\quad 
\lim_{{i_k} \to \infty}\frac{{q'_{i_k}}}{q_{i_k} + 2q'_{i_k}}= \frac{\lambda - 1}{\lambda}.
\end{equation}
Analogously in the second case.

By lemma \ref{cases} we can find sequence $\{i_k\}_{k = 1}^{\infty}$ such that for all
$i_k$ one of the following cases is realized:
\begin{enumerate}
\item Case 1 of Lemma holds, $x < \frac{p_{i_k} + p'_{i_k}}{q_{i_k} + q'_{i_k}}$.
\item Case 1 of Lemma holds, $x > \frac{p_{i_k} + p'_{i_k}}{q_{i_k} + q'_{i_k}}$
\item Case 2 of Lemma holds, $x < \frac{p_{i_k} + p'_{i_k}}{q_{i_k} + q'_{i_k}}$
\item Case 2 of Lemma holds, $x > \frac{p_{i_k} + p'_{i_k}}{q_{i_k} + q'_{i_k}}$
\end{enumerate}

Let us consider for example fist of these cases. As case 1 of lemma \ref{cases} holds,
then one of equalities (\ref{e1}), (\ref{e2}) is satisfied. As
$x < \frac{p_{i_k} + p'_{i_k}}{q_{i_k} + q'_{i_k}}$, then 
one of equalities (\ref{e5}), (\ref{e6}) is satisfied.
In such way we get four variants:
\begin{itemize}
\item
\begin{equation}\label{e9}
\lim_{{i_k} \to \infty}\frac{{q_{i_k}}'}{q_{i_k} + q'_{i_k}} = \frac{\lambda - 1}{\lambda},
\quad \lim_{{i_k} \to \infty}\frac{{q_{i_k}}}{2q_{i_k} + q'_{i_k}} = \frac{\lambda - 1}{\lambda}.
\end{equation}
Since
$$
\frac{2q_{i_k} + q'_{i_k}}{{q_{i_k}}} = 1 + \frac{q_{i_k} + q'_{i_k}}{{q_{i_k}}'}
$$
we get from (\ref{e9}) incorrect equality 
$$
\frac{\lambda}{\lambda - 1} = 1 + \frac{\lambda}{\lambda - 1}.
$$
\item
\begin{equation}\label{e10}
\lim_{{i_k} \to \infty}\frac{{q_{i_k}}'}{q_{i_k} + q'_{i_k}} = \frac{1}{\lambda},
\quad \lim_{{i_k} \to \infty}\frac{{q_{i_k}}}{2q_{i_k} + q'_{i_k}} = \frac{\lambda - 1}{\lambda}.
\end{equation}
Analogously to the previous case we get from (\ref{e10}) following equation:
$$
\lambda^2 - \lambda - 1 = 0.
$$
But $\lambda$ is a root of irreducible equation of degree $3$.
\item
\begin{equation}\label{e11}
\lim_{{i_k} \to \infty}\frac{{q_{i_k}}'}{q_{i_k} + q'_{i_k}} = \frac{\lambda - 1}{\lambda},
\quad \lim_{{i_k} \to \infty}\frac{{q_{i_k}}}{2q_{i_k} + q'_{i_k}} = \frac{1}{\lambda}.
\end{equation}
From (\ref{e11}) we get following equation:
$$
\lambda^2 - 3 \lambda +1 = 0.
$$ 
But $\lambda$ is a root of irreducible equation of degree $3$.
\item
\begin{equation}\label{e12}
\lim_{{i_k} \to \infty}\frac{{q_{i_k}}'}{q_{i_k} + q'_{i_k}} = \frac{1}{\lambda},
\quad \lim_{{i_k} \to \infty}\frac{{q_{i_k}}}{2q_{i_k} + q'_{i_k}} = \frac{1}{\lambda}.
\end{equation}
From (\ref{e12}) we get incorrect equation:
$$
\lambda = 1 + \lambda.
$$ 
\end{itemize} 
So we get contradiction in all of examined cases. The rest of the cases can be examined analogously.
\hfill $\blacksquare$

\end{document}